\documentclass[11pt]{article}
\usepackage[a4paper]{anysize}\marginsize{2cm}{2cm}{2cm}{2cm}
\pdfpagewidth=\paperwidth \pdfpageheight=\paperheight
\usepackage{amsfonts,amssymb,amsthm,amsmath,eucal}
\usepackage{bbm}
\usepackage{tikz} 
\usepackage{mathrsfs}

\theoremstyle{plain}
\newtheorem{thm}{Theorem}[section]
\newtheorem{theorem}[thm]{Theorem}

\newtheorem{conjecture}[thm]{Conjecture}

\theoremstyle{definition}
\newtheorem{definition}[thm]{Definition}

\newtheorem{example}[thm]{Example}

\newtheorem{problem}[thm]{Problem}

\newtheorem{thevarthm}[thm]{\varthmname}

\newenvironment{varthm*}[1]{\trivlist\item[]{\bf #1.}\it}{\endtrivlist}

\renewcommand\geq{\geqslant}

\renewcommand\leq{\leqslant}

\newcommand\be{\begin{eqnarray*}}
\newcommand\ee{\end{eqnarray*}}

\newcommand\R{\mathbb R}
\newcommand\C{\mathbb C}

\newcommand\K{\mathbb K}
\renewcommand\P{\mathbb P}

\newcommand\newop[2]{\def#1{\mathop{\rm #2}\nolimits}}
\newop\log{log}
\newop\ord{ord}
\newop\Gal{Gal}
\newop\SL{SL}
\newop\Bl{Bl}
\newop\mult{mult}
\newop\imult{imult}
\newop\mass{mass}
\newop\Ass{Ass}
\newop\div{div}
\newop\codim{codim}
\newop\sing{sing}
\newop\Zeroes{Zeroes}

\newcommand\alphahat{\widehat{\alpha}}
\newcommand\eps{\varepsilon}

\def\keywordname{{\bfseries Keywords}}%
\def\keywords#1{\par\addvspace\medskipamount{\rightskip=0pt plus1cm
\def\and{\ifhmode\unskip\nobreak\fi\ $\cdot$
}\noindent\keywordname\enspace\ignorespaces#1\par}}
\def\subclassname{{\bfseries Mathematics Subject Classification
(2000)}\enspace}
\def\subclass#1{\par\addvspace\medskipamount{\rightskip=0pt plus1cm
\def\and{\ifhmode\unskip\nobreak\fi\ $\cdot$
}\noindent\subclassname\ignorespaces#1\par}}

\newcommand\beginproof[1]{\trivlist\item[\hskip\labelsep{\em #1.}]}
\newcommand\proofof[1]{\beginproof{Proof of #1}}

\def\endproof{\hspace*{\fill}\endproofsymbol\endtrivlist}

\def\endproofsymbol{\frame{\rule[0pt]{0pt}{6pt}\rule[0pt]{6pt}{0pt}}}

\begin{document}

\title{Waldschmidt constants for Stanley-Reisner ideals of a class of graphs}
\author{Tomasz Szemberg and Justyna Szpond}
\maketitle

\abstract{In the present note we study Waldschmidt constants of Stanley-Reisner ideals of a hypergraph
   and a graph with vertices forming a bipyramid over a planar $n$--gon. The case of the hypergraph
   has been studied by Bocci and Franci. We reprove their main result. The case of the graph is new.
   Interestingly, both cases provide series of ideals with Waldschmidt constants descending to $1$.
   It would be interesting to known if there are bounded ascending sequences of Waldschmidt constants.}

\section{Introduction}
\label{SSsec:1}
   The following problem has attracted considerable attention in commutative
   algebra and algebraic geometry in the past two decades.
   \begin{problem}[Containment problem]\label{pro:containment}
      Let $I$ be a homogeneous ideal in the polynomial ring $\K[x_0,\ldots,x_N]$,
      where $\K$ is a field. Decide for which integers $m$ and $r$ there is the
      containment
      \begin{equation}\label{eq:containment}
         I^{(m)}\subset I^r
      \end{equation}
      between the symbolic and ordinary powers of the ideal $I$.
   \end{problem}
   We recall that for $m\geq 0$ the $m$-th symbolic power of $I$ is defined as
   \begin{equation}\label{eq:symbolic power}
      I^{(m)}=\bigcap_{P\in\Ass(I)}\left(I^mR_P\cap R\right),
   \end{equation}
   where $\Ass(I)$ is the set of associated primes of $I$.
   At the beginning of the Millennium, Ein, Lazarsfeld and Smith
   in characteristic zero \cite{ELS01} and Hochster and Huneke
   in positive characteristic \cite{HoHu02} proved striking uniform
   answers to Problem \ref{pro:containment} to the effect that
   the containment in \eqref{eq:containment} holds for all
   \begin{equation}\label{eq:ELS HH}
      m\geq hr,
   \end{equation}
   where $h$ is the maximum of heights of all associated primes of $I$.
   In geometric terms it means that $h$ is the codimension of the smallest embedded
   component of the set $\Zeroes(I)$. In particular, the containment in \eqref{eq:containment}
   holds for all $I$ with $m\geq Nr$.

   It is natural to wonder to what extend the particular bound in \eqref{eq:ELS HH}
   is sharp. In order to study this question Bocci and Harbourne introduced in \cite{BocHarJAG10}
   a number of \emph{asymptotic invariants} attached to $I$. In the present note
   we focus on one of them. Let $\alpha(I)$ denote the smallest degree of a non-zero
   element in $I$, this is the \emph{initial degree} of $I$. Then, the \emph{Waldschmidt constant}
   of $I$ is asymptotically defined as
   \begin{equation}\label{eq:Waldschmidt constant}
      \alphahat(I)=\lim_{m\to\infty}\frac{\alpha(I^{(m)})}{m}.
   \end{equation}
   It is well known, see e.g. proof of \cite[Lemma 2.3.1]{BocHarJAG10},
   that $\alphahat(I)=\inf_{m\geq 1}\frac{\alpha(I^{(m)})}{m}$.

   Interestingly, Waldschmidt constants were introduced long before
   the Problem \ref{pro:containment} has been considered in the realms
   of complex analysis, see \cite{Wal77} and our note \cite{MSS17}
   for recent account. These invariants are very hard to compute in
   general. In fact, a number of important conjectures can be
   expressed in terms of Waldschmidt constants. By the way of an example
   we mention here only the following one.
\begin{conjecture}[Nagata]
   Let $I$ be the ideal defining $s\geq 10$ \emph{very general} points
   in $\P^2$. Then
   $$\alphahat(I)=\sqrt{s}.$$
\end{conjecture}
   Our research here has been motivated by \cite{BocFra16}, where
   Bocci and Franci initiated the study of Waldschmidt constants
   of monomial ideals determined by some combinatorial data.
   They have computed Waldschmidt constants of Stanley-Reisner
   ideals of bipyramids (see Section \ref{SSsubs:2.2}).
   We provide here an alternative, more elementary proof of their result and
   study Stanley-Reisner ideals of graphs determined by vertices
   of bipyramides. Our main result is Theorem \ref{thm:Waldschmidt for D_n}.

\section{Bipyramids revisited}
\label{SSsec:2}
   We begin by recalling briefly basic notions from combinatorial
   algebra relevant in this note, for more detailed account see
   the very nice surveys \cite{FHM13} and \cite{FMS14}. The Stanley-Reisner ideals
   introduced here have traditionally provided a rich source
   of non-trivial examples.
\subsection{Stanley-Reisner theory}
\label{SSsubs:2.1}
\begin{definition}[Simplicial complex]
   A \emph{simplicial complex} $\Delta$ on a finite set $V$
   is a collection of subsets $\sigma\subset V$ such that
   the containment $\sigma\in\Delta$ implies $\tau\in\Delta$
   for all subsets $\tau\subset\sigma$.
\end{definition}
   For the set $V=\left\{0,1,\ldots,N\right\}$, we can
   naturally identify any subset $\sigma\subset V$ with
   a squarefree monomial
   $$x_{\sigma}=\prod_{i\in\sigma}x_i\in\K[x_0,\ldots,x_N].$$
   The key observation of the Stanley-Reisner theory is that there
   is a bijective correspondence between squarefree monomial ideals
   and simplicial complexes.
\begin{definition}[Stanley-Reisner ideal]
   The \emph{Stanley-Reisner ideal} of the simplicial complex $\Delta$
   is the monomial ideal
   $$I_{\Delta}=\langle x_{\tau}:\; \tau\notin\Delta\rangle.$$
\end{definition}
   There is a big advantage of working with symbolic powers of monomial ideals
   rather than symbolic powers of arbitrary ideals. It follows from the
   following extremely useful result that one can avoid localizations, see \cite[Theorem 3.7]{CEHH16}
   and \cite[Theorem 2.5]{BCGHJNSvTV16}.
\begin{theorem}[Symbolic powers of monomial ideals]\label{thm:symbolic powers of monomial}
   Let $I\subset \K[x_0,\ldots,x_N]$ be a monomial ideal with minimal primary decomposition
   $$I=P_1\cap\ldots\cap P_s.$$
   Then, for all $m\geq 1$ there is
   $$I^{(m)}=P_1^m\cap\ldots\cap P_s^m.$$
\end{theorem}
\subsection{Bipyramids}
\label{SSsubs:2.2}
   Following Bocci and Franci \cite{BocFra16}, for $n\geq 3$, we define
   a \emph{bipyramide} $B_n$ over an $n$--gon $\Gamma_n$ as the convex hull of the following set of points
   $$\left\{(0,1),(1,0),(\eps,0),(\eps^2,0),\ldots,(\eps^{n-1},0),(0,-1)\right\}\subset \C\times\R,$$
   where $\eps$ is a primitive root of $1$ of order $n$ and $\Gamma_n$ has vertices in the plane $y=0$.
   Thus a bipyramid is a polytop. Numbering its vertices as follows
   $$P_0=(0,1), P_k= (\eps^k,0)\;\mbox{ for }\;k=1,\ldots, n,\; P_N=P_{n+1}=(0,-1)$$
   and assigning to each face of $B_n$ the set of its vertices, we obtain
   a simplicial complex with $V=\left\{0,1,\ldots,N\right\}$.
   Thus its Stanley-Reisner ideal is
   $$I_{B_n}=\langle x_0x_N, x_ix_j \mbox{ with }\;1\leq i<j\leq n\; \mbox{ and } P_iP_j \mbox{ not an edge of }\Gamma_n\rangle.$$
   For $i=1,\ldots, n$, let
   \begin{equation}\label{eq:notation T_i}
      S_i=\left\{x_i,x_{i+1},\ldots,x_{i+n-2}\right\}\; \mbox{ and }\; T_i=\left\{x_i,x_{i+1},\ldots,x_{i+n-3}\right\},
   \end{equation}
   where the indices are numbered so that $x_{n+i}=x_i$ for $i=1,\ldots,N$.
   It is easy to check that
   \begin{equation}\label{eq:I_B_n primary decomposition}
      I_{B_n}=\bigcap_{i=1}^{n} \langle x_0,T_i\rangle
         \cap\bigcap_{i=1}^{n} \langle x_N,T_i\rangle
   \end{equation}
   is the primary decomposition of $I_{B_n}$, see \cite[Proposition 3.2]{BocFra16}.
\begin{example}
   For $n=4$ the bipyramid $B_4$ is indicated in Figure \ref{fig:B4}. We have
   \begin{equation*}
      I_{B_4}=\langle x_0x_5,x_1x_3,x_2x_4 \rangle
   \end{equation*}
   \begin{equation*}
   \begin{split}
   I_{B_4} & =\langle x_0,x_1,x_2 \rangle \cap \langle x_0,x_2,x_3 \rangle \cap \langle x_0,x_3,x_4 \rangle \cap \langle x_0,x_1,x_4 \rangle \cap \\ &\;\, \cap
      \langle x_1,x_2,x_5 \rangle \cap \langle x_2,x_3,x_5 \rangle \cap \langle x_3,x_4,x_5 \rangle \cap \langle x_1,x_4,x_5 \rangle.
   \end{split}
   \end{equation*}
\begin{figure}[h]
\begin{center}
\includegraphics[scale=.65]{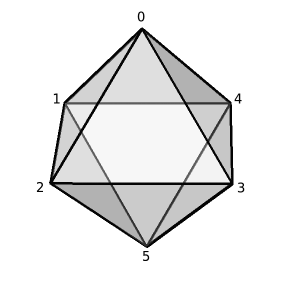}
\end{center}
\caption{The bipyramid $B_4$}
\label{fig:B4}
\end{figure}
\end{example}
   The main result of \cite{BocFra16} is the following Theorem (\cite[Theorem 1.1]{BocFra16}).
\begin{theorem}\label{thm:BocFra}
   For any $n\geq 4$, the Waldschmidt constant of the Stanley-Reisner ideal $I_{B_n}$
   of a bipyramid $B_n$ is
   $$\alphahat(I_{B_n})=\frac{n}{n-2}.$$
   For $n=3$ there is $I_{B_3}=\langle x_0x_4 \rangle$ and hence $\alphahat(I_{B_3})=2$.
\end{theorem}
   This Theorem has been already reproved in \cite[Theorem 6.10]{BCGHJNSvTV16},
   the authors appeal however to fractional chromatic numbers of hypergraphs
   and use the advanced machinery developed in their paper. We provide here,
   as an alternative, yet another, fairly elementary proof.
\proofof{Theorem \ref{thm:BocFra}}
   By definition \eqref{eq:Waldschmidt constant}, the Waldschmidt constant is a limit
   of a sequence, hence it can be computed by an arbitrary subsequence. We use
   the subsequence indexed by $s(n-2)$ for $s\geq 1$.\\
   Our first observation is that
   $$(x_1\cdot\ldots\cdot x_n)\in I^{(n-2)}.$$
   Indeed, combining Theorem \ref{thm:symbolic powers of monomial} and \eqref{eq:I_B_n primary decomposition},
   we see that $I^{(n-2)}$ is the intersection of ideals of the type
   $$\langle x_u,T_i \rangle^{(n-2)},$$
   where $u\in\left\{0,N\right\}$. Thus, clearly
   $$x_i\cdot x_{i+1}\cdot\ldots\cdot x_{i+n-3}\;| x_1\cdot\ldots\cdot x_n.$$
   Since $\deg(x_1\cdot\ldots\cdot x_n)=n$, we have
   \begin{equation}\label{eq:upper bound}
      \alphahat(I_{B_n})\leq \frac{n}{n-2}.
   \end{equation}
   Turning to the reverse inequality, assume by the way of contradiction
   that there is a monomial $f=x_0^{a_0}\cdot x_1^{a_1}\cdot\ldots\cdot x_{n+1}^{a_{n+1}}$
   of degree $\leq sn-1$ in $I^{(s(n-2))}$, i.e.
   \begin{equation}\label{eq:zzz-1}
   \sum_{i=0}^N a_i\leq sn-1.
   \end{equation}
   Since $f$ is contained in all ideals $\langle x_u,T_i \rangle^{s(n-2)}$
   with $u\in\left\{0,N=n+1\right\}$ and $i=1,\ldots,n$, we obtain
   $2n$ inequalities of the type
   \begin{equation}\label{eq:zzz}
      a_u+a_i+a_{i+1}+\ldots +a_{i+n-3}\geq s(n-2).
   \end{equation}
   Summing these inequalities, we get
   $$n(a_0+a_N)+2(n-2)\sum_{i=1}^n a_i\geq 2ns(n-2).$$
   Since $n\geq 4$, the left hand side is bounded from above by $2(n-2)\sum_{i=0}^N a_i$.
   Taking \eqref{eq:zzz-1} into account, we obtain
   $$2ns(n-2)\leq 2(n-2)(sn-1),$$
   which is a clear contradiction. Hence we conclude that
   \begin{equation}\label{eq:lower bound}
      \frac{\alpha(I^{(s(n-2))})}{s(n-2)}\geq \frac{sn}{s(n-2)}=\frac{n}{n-2}\;\mbox{ for all }\; s\geq 1.
   \end{equation}
   Combining \eqref{eq:upper bound} with \eqref{eq:lower bound} we obtain the assertion.
\endproof
\section{Bipyramidal graph}
\label{SSsec:3}
   In this section we consider a graph, rather than a hypergraph, defined by vertices
   of a bipyramid. To be more precise, we define the bipyramidal graph $D_n$ as the set of vertices
   $V=\left\{P_0,P_1,\ldots,P_n,P_N\right\}$ with $N=n+1$ together with the set of edges
   $$E=\left\{ P_0P_i,\; P_NP_i,\; P_1P_2,\; P_2P_3,\;\ldots,\;P_{n-1}P_n,\;P_nP_1,\;P_0P_N\mbox{ for }\; i=1,\ldots,n\right\}.$$
\begin{example}
   For $n=4$ the graph is indicated in Figure \ref{fig:D4}. We have
   \begin{equation*}
      I_{D_4}=\langle x_1x_3,x_2x_4,x_0x_1x_2,x_0x_2x_3,x_0x_3x_4,x_0x_4x_1,x_1x_2x_5,x_2x_3x_5,x_3x_4x_5,x_4x_1x_5 \rangle
   \end{equation*}
   \begin{equation*}
   \begin{split}
   I_{D_4} & =\langle x_0,x_5,x_1,x_2 \rangle \cap \langle x_0,x_5,x_2,x_3 \rangle \cap \langle x_0,x_5,x_3,x_4 \rangle \cap \langle x_0,x_5,x_1,x_4 \rangle \cap \\ &\;\, \cap
      \langle x_1,x_2,x_3,x_5 \rangle \cap \langle x_2,x_3,x_4,x_5 \rangle \cap \langle x_3,x_4,x_1,x_5 \rangle \cap \langle x_1,x_2,x_4,x_5 \rangle \cap \\ &\;\, \cap
      \langle x_0,x_1,x_2,x_3 \rangle \cap \langle x_0,x_2,x_3,x_4 \rangle \cap \langle x_0,x_3,x_4,x_1 \rangle \cap \langle x_0,x_1,x_2,x_4 \rangle \cap \\ &\;\, \cap
      \langle x_1,x_2,x_3,x_4 \rangle.
   \end{split}
   \end{equation*}
\begin{figure}[h]
\begin{center}
\includegraphics[scale=.65]{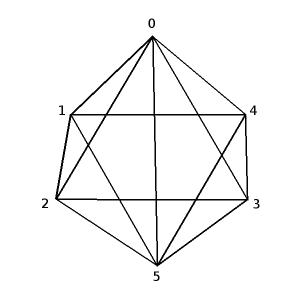}
\end{center}
\caption{The bipyramidal graph $D_4$}
\label{fig:D4}
\end{figure}
\end{example}
\begin{theorem}\label{thm:Waldschmidt for D_n}
   For the Stanley-Reisner ideal $I_{D_n}$ of the bipyramidal graph $D_n$ we have
   $$\alphahat(I_{D_n})=\frac{n+2}{n}.$$
\end{theorem}
\proof
   Note to begin with that using the notation in \eqref{eq:notation T_i}
   \begin{equation}\label{eq:D_n primary decomposition}
   \begin{split}
   I_{D_n}=\bigcap_{i=1}^{n}\langle x_0, x_N, T_i\rangle \cap
           \bigcap_{i=1}^{n}\langle x_0, S_i\rangle \cap
           \bigcap_{i=1}^{n}\langle x_N, S_i\rangle \cap
           \langle x_1,\ldots,x_n \rangle
   \end{split}
   \end{equation}
   is the primary decomposition.\\
   It follows that $x_0\cdot x_N\cdot x_1\cdot\ldots\cdot x_n\in I_{D_n}^{(n)}$,
   hence
   \begin{equation}\label{eq:upper bound D_n}
      \alphahat(I_{D_n})\leq \frac{n+2}{n}.
   \end{equation}
   Turning to the lower bound, we study the sequence of symbolic powers
   of $I_{D_n}$ indexed by multiples $sn$ of $n$ for $s\geq 1$.
   We assume that there is an $s$ such that $I_{D_n}^{(ns)}$ contains
   a monomial $g=x_0^{a_0}\cdot x_N^{a_N}\cdot x_1^{a_1}\cdot\ldots\cdot x_n^{a_n}$ of degree
   \begin{equation}\label{eq:D1}
      a_0+a_N+a_1+\ldots+a_n\leq s(n+2)-1.
   \end{equation}
   Since $g\in\langle x_1,\ldots,x_n\rangle$, we also have
   \begin{equation}\label{eq:D2}
      a_1+\ldots+a_n\geq sn.
   \end{equation}
   It follows from \eqref{eq:D1} and \eqref{eq:D2} that
   \begin{equation}\label{eq:D3}
      a_0+a_{N}\leq 2s-1.
   \end{equation}
   Since there is in the decomposition \eqref{eq:D_n primary decomposition}
   an ideal which misses arbitrary two consecutive indices in the set $\left\{0,1,\ldots,n,N=(n+1)\right\}$,
   we obtain by the same token that
   \begin{equation}\label{eq:D4}
      a_i+a_{i+1}\leq 2s-1.
   \end{equation}
   for $i=0,\ldots,n$.
   Summing up altogether $(n+2)$ inequalities in \eqref{eq:D3} and \eqref{eq:D4},
   we obtain
   \begin{equation}\label{eq:D5}
      a_0+a_N+a_1+\ldots+a_n\leq (n+2)(s-\frac12).
   \end{equation}
   On the other hand, since $g$ is an element in all ideals in the decomposition \eqref{eq:D_n primary decomposition},
   we obtain, analogously to \eqref{eq:D2}
   \begin{equation}\label{eq:D6}
      a_i+\ldots+a_{n+i}\geq sn.
   \end{equation}
   for all $i=0,\ldots,N$, of course with the convention that the indices
   are taken modulo $(n+2)$.
   Summing up these inequalities we get
   \begin{equation}\label{eq:D7}
      a_0+a_N+a_1+\ldots+a_n\geq (n+2)s.
   \end{equation}
   Inequalities \eqref{eq:D5} and \eqref{eq:D7} give the desired contradiction,
   implying that all polynomials in $I_{D_n}^{ns}$ have degree al least $s(n+2)$.
   This, in turn, implies that
   \begin{equation}\label{eq:lower bound D_n}
      \alphahat(I_{D_n})=\lim_{s\to\infty}\frac{\alpha(I_{D_n}^{(ns)})}{ns}\geq \frac{s(n+2)}{sn}=\frac{n+2}{n}.
   \end{equation}
   Thus \eqref{eq:upper bound D_n} and \eqref{eq:lower bound D_n} establish the assertion and we are done.
\endproof

\paragraph{Acknowledgements.}
   Our research was partially supported by National Science Centre, Poland, grant
   2014/15/B/ST1/02197.


\bigskip \small

\bigskip
   Tomasz Szemberg,
   Department of Mathematics, Pedagogical University of Cracow,
   Podchor\c a\.zych 2,
   PL-30-084 Krak\'ow, Poland

\nopagebreak
   \textit{E-mail address:} \texttt{tomasz.szemberg@gmail.com}

\bigskip
   Justyna Szpond,
   Department of Mathematics, Pedagogical University of Cracow,
   Podchor\c a\.zych 2,
   PL-30-084 Krak\'ow, Poland

\nopagebreak
   \textit{E-mail address:} \texttt{szpond@up.krakow.pl}

\end{document}